\documentclass[12pt]{article}
\usepackage{amsmath,amssymb,amsthm}
\usepackage{graphics,epsfig,calc}
\usepackage{amsmath}
\usepackage{latexsym,epsfig,bm,amssymb}
\usepackage{color}
\usepackage{amsthm,mathrsfs}

\usepackage{mathptmx}

\newcommand{\be}{\begin{equation}}
\newcommand{\ee}{\end{equation}}
\newcommand{\beqn}{\begin{eqnarray}}
\newcommand{\eeqn}{\end{eqnarray}}

\newcommand{\ba}{\begin{array}}
\newcommand{\ea}{\end{array}}
\newcommand{\bo}{{\hfill\loota}}

\newcommand{\loota}{\hbox{\enspace{\vrule height 7pt depth 0pt width 7pt}}}

\newcommand{\cE}{{\cal E}}

\newcommand{\cH}{{\cal H}}

\newcommand{\cV}{{\cal V}}

\newcommand{\cX}{{\cal X}}

\newcommand{\ds}{\displaystyle}
\newcommand{\fr}{\frac}

\newcommand{\de}{\delta}\newcommand{\De}{\Delta}

\newcommand{\Lam}{\Lambda}
\newcommand{\lam}{\lambda}

\newcommand{\om}{{\omega}}
\newcommand{\si}{\sigma}
\newcommand{\vp}{\varphi}

\newcommand{\la}{\label}

\newcommand{\ov}{\overline}

\newcommand{\re}{\ref}

\newcommand{\ti}{\tilde}

\newcommand\C{{\mathbb C}}
\newcommand\R{{\mathbb R}}

\newcommand{\5}{{\hspace{0.5mm}}}

\newcommand{\sgn}{\mathop{\rm sgn}\nolimits}

\newcommand{\rRe}{{\rm Re\5}}
\newcommand{\rIm}{{\rm Im\5}}

\renewcommand{\Pr}{\hspace{-6mm}{\bf Proof~}}

\newcommand{\Ker}{{\rm Ker\5}}
\newcommand{\Ran}{{\rm Ran\5}}

\renewcommand{\theequation}{\thesection.\arabic{equation}}
\newtheorem{theorem}{Theorem}[section]
\renewcommand{\thetheorem}{\arabic{section}.\arabic{theorem}}
\newtheorem{definition}[theorem]{Definition}

\newtheorem{lemma}[theorem]{Lemma}
\newtheorem{example}[theorem]{Example}
\newtheorem{remark}[theorem]{Remark}
\newtheorem{remarks}[theorem]{Remarks}
\newtheorem{cor}[theorem]{Corollary}
\newtheorem{pro}[theorem]{Proposition}

\newcommand{\bd}{\begin{definition}}
 \newcommand{\ed}{\end{definition}}
\newcommand{\bt}{\begin{theorem}}
 \newcommand{\et}{\end{theorem}}
\newcommand{\bqt}{\begin{qtheorem}}
 \newcommand{\eqt}{\end{qtheorem}}

\newcommand{\bp}{\begin{pro}}
 \newcommand{\ep}{\end{pro}}

\newcommand{\bl}{\begin{lemma}}
 \newcommand{\el}{\end{lemma}}
\newcommand{\bc}{\begin{cor}}
 \newcommand{\ec}{\end{cor}}

\newcommand{\bex}{\begin{example}}
 \newcommand{\eex}{\end{example}}
\newcommand{\bexs}{\begin{examples}}
 \newcommand{\eexs}{\end{examples}}

\newcommand{\bexe}{\begin{exercice}}
 \newcommand{\eexe}{\end{exercice}}

\newcommand{\br}{\begin{remark} }
 \newcommand{\er}{\end{remark}}
\newcommand{\brs}{\begin{remarks}}
 \newcommand{\ers}{\end{remarks}}

\begin{document}
\begin{titlepage}
\begin{center}
{\Large\bf
On eigenfunction expansion of solutions
\bigskip\\
to the Hamilton equations}
\end{center}
\vspace{2cm}
 \begin{center}
{\large A. Komech}
\footnote{Supported partly
by Alexander von Humboldt Research Award,
Austrian Science Fund (FWF): P22198-N13,
and the grant of the Russian Foundation for Basic Research.}
\\
{\it Faculty of Mathematics of Vienna University\\
and Institute for Information Transmission Problems RAS} \\
e-mail:~alexander.komech@univie.ac.at
\bigskip\\
{\large E. Kopylova}
\footnote{
Supported partly by
Austrian Science Fund (FWF): M1329-N13,
and the grant of the Russian Foundation for Basic Research.}
\\
{\it Faculty of Mathematics of Vienna University\\
and Institute for Information Transmission Problems RAS} \\
 e-mail:~elena.kopylova@univie.ac.at
\end{center}
\vspace{1cm}

 \begin{abstract}

We establish a spectral representation for solutions to
linear
Hamilton equations 
with positive definite energy
in a~Hilbert space.
Our approach is a special version of M.~Krein's 
spectral theory of $J$-selfadjoint operators in 
the Hilbert spaces with indefinite metric.
Our main result is an application to
the eigenfunction expansion for
the linearized relativistic Ginzburg--\allowbreak Landau equation.

\bigskip\bigskip\bigskip


\end{abstract}

\end{titlepage}

\section{Introduction}

We consider complex linear Hamilton operators in a complex Hilbert space $\cX$,
\begin{equation}\la{Hs}
A=JB~,~~~ \mbox{ where }~~~B^*=B~,~~~~J^*=-J~,~~J^2=-1~.
\end{equation} 
In other words, $A$ is a $J$-selfadjoint operator \cite{L1981}.
The  selfadjoint operator
 $B$ is defined on a dense domain $D(B)\subset\cX$.
Our aim is a spectral representation for solutions to
 the equation
\begin{equation}\la{NHl}
\dot X(t)=AX(t)~,~~~~~~~~~~~~t\in\R~.
\end{equation}
Our main goal is an application to
the eigenfunction expansion for
the linearized relativistic Ginzburg--\allowbreak Landau equation.

In the simplest case, when $J=i$,  the
solutions are given by 
$\ds X(t)=e^{iBt}X(0)$.
A more general  `commutative case', when $JB=BJ$,
reduces to the case $J=i$ since $JB=iB_1$, where $B_1=-iJB$ is the selfadjoint operator.
However, 
$JB\ne BJ$ for linearizations of
$U(1)$-invariant
nonlinear Schr\"odinger equations as we show in Appendix B.

The complex Hamilton operators arise as the complexification of
a linearization of nonlinear Hamilton equations
in a real Hilbert space $\cX_r$, 
\begin{equation}\la{NH}
\dot \psi_r(t)=J_r D\cH(\psi_r(t))~,
\end{equation}
where
$\cH$ is the Hamilton functional,
$J_r^*=-J_r$, $J_r^2=-1$,
and
$D\cH$
stands for the differential which is a linear operator
over real numbers.
The linearization of (\re{NH}) at a stationary state $S\in \cX_r$
reads as
\begin{equation}\la{NHlr}
\dot X_r(t)=A_rX_r(t)~,\qquad A_r=J_rB_r~, 
\end{equation}
where $B_r:=DD\cH(S)$
is a real symmetric operator. This is the Hamilton system with the
Hamilton functional
\begin{equation}\la{hamr}
\cH_0(X_r)=\fr 12\langle B_rX_r,X_r\rangle_r~,
\end{equation}
where the brackets denote the scalar product on $\cX_r$.
The complexification of the operator
$A_r$ in the complex space $\cX=\cX_r\oplus i\cX_r$
reads as (\re{Hs}), where $J$, $B$ are the complexifications of
$J_r$, $B_r$ respectively.
Moreover, $B$ is complex selfadjoint on $\cX$ if $B_r$ is real 
selfadjoint on~$\cX_r$.
\medskip

The representation of solutions
to the Hamilton equations
in the form of oscillatory integrals
is indispensable in the proof of the dispersion decay
for
linearized equations
in the theory of asymptotic stability of solitary waves
of nonlinear Schr\"odinger, wave, Klein--Gordon, Maxwell, Dirac
and relativistic Ginzburg--\allowbreak Landau
equations, which were extensively developed the last two decades
\cite{BC11}--\cite{BS03}, \cite{C01}--\cite{CT09},
\cite{IKS12}--\cite{KKm11}, \cite{MM08}--\cite{PW94},
and \cite{Sigal93}--\cite{T02}.
However, many features of 
these representations were not justified up to now,
since the generators of the 
linearized equations can be non symmetric,
as noted  in Introduction of \cite{C11H}. 
In particular, the eigenfunction expansion, used in 
\cite{KK11}, was not justified with detail.
\medskip

We fill this gap
in the simplest 
case of positive definite `energy operators' 
$B$ satisfying
the following spectral condition:
\beqn\la{P}
\mbox{\rm\bf Condition ~~~P$\bf_+$}~~~~~~~~~~~~~~~~~~~~~~~~~~~~~~~
\sigma (B)\subset [\delta,\infty)~,
~~~~\delta>0~.~~~~~~~~~~~~~~~~~
\eeqn
Equivalently, $B$ is invertible in $\cX$ and $B> 0$.
This framework is sufficient for our main application to
the linearized relativistic Ginzburg--\allowbreak Landau equation.
\medskip\\
Our results are the following:
\medskip\\
$\bullet$
The similarity of $iA$ to a selfadjoint operator.
\medskip\\
$\bullet$
The existence and uniqueness of generalized solutions
to (\re{NHl})
for all initial
states $X$ with finite
energy $\langle B X,X\rangle$, where $\langle \cdot,\cdot \rangle$ stands for
the scalar product on $\cX$.
\medskip\\
$\bullet$
A spectral representation of solutions
for complex and real Hamilton generators.
\medskip\\
$\bullet$
Our main result is 
 the eigenfunction expansion
\begin{equation}\la{uns30i}
\left(
\ba{c}
\psi(t)\\
\dot \psi(t)
\ea
\right)
=
\sum_{-N}^N e^{-i\omega_k t} C_k a_k+
\int_{|\omega|\ge m} e^{-i\omega t} C(\omega)a_\omega~d\omega
\end{equation}
for solutions to the
linearized
relativistic nonlinear Ginz\-burg--\allowbreak Landau
equation \cite{KK11}.
Here 
$a_k$ are the eigenfunctions of the corresponding Klein-Gordon
generator
\begin{equation}\la{AB0i}
A=\left(
\ba{cc}
0   & 1 \\
-H_0 & 0
\ea
\right)~,
\end{equation}
where $H_0:=-\fr{d^2}{dx^2}+m^2+ V_0(x)$. Finally, 
$a_\omega$ are the 
generalized eigenfunctions of the continuous spectrum of $A$.
\medskip\\
Such eigenfunction expansions were used in 
 \cite{BP95,BS03,KK11}
for the calculation of `Fermi Golden Rule' (FGR) in the context of 
nonlinear 
Schr\"odinger and Klein-Gordon equations.
This is a nondegeneracy condition
 introduced
in \cite{Sigal93}
in the framework of nonlinear wave and
Schr\"odinger equations.
The condition
provides a strong coupling of discrete
and continuous spectral components of solutions, 
which provides the
energy radiation to infinity and results in the
 asymptotic stability of solitary
waves.
The calculation of FGR, as given in \cite{BP95,BS03,KK11}, relies
on eigenfunction expansions of type 
(\re{uns30i}).
Our main Theorem \re{tmain} justifies the eigenfunction expansion
\cite[(5.14)]{KK11}
which was not proved  with detail. This justification
was
one of our main motivation
in writing the present paper.
\medskip

Let us comment on our approach.
First, we 
reduce the problem to a selfadjoint generator
defined uniquely by $A$ justifying the classical  M. Krein transformation
under condition
(\re{P}). 
This reduction
is a special version of spectral theory
of $J$-selfadjoint operators in the Hilbert spaces with an indefinite metric
 \cite{AI1989,KL1963}. 
We plan to extend elsewhere
these methods and results to more general degenerated case when $\Ker B\ne 0$.

Second, we apply the abstract spectral theory to the operator
(\re{AB0i}) and develop
a modification of the eigenfunction expansion theory \cite[pp 114-115]{RS3}
for the reduced selfadjoint operator,
see Remark \re{rver}.
Finally, we apply the theory of PDO
to
deduce (\re{uns30i}) from this
modification.

One of the  novelties here is a vector-valued 
treatment 
 of the convergence
of the integral over the continuous spectrum in (\re{uns30i}).
Namely, we show 
that the integral is the limit of the 
corresponding integrals over $m\le|\omega|\le M$   as $M\to\infty$
in the Sobolev space 
$H^1(\R)\oplus L^2(\R)$.
In its own turn, the integral over $m\le|\omega|\le M$ is absolutely 
converging in the weighted $L^2$-space with the weight $(1+|x|)^{-s}$
where $s>1$.

\medskip

Let us comment on related works. 
Some spectral properties of the Hamilton non-selfadjoint 
operators were studied by 
V. Buslaev and G. Perelman \cite{BP93,BP95,BS03},
M.B. Erdogan and W. Schlag \cite{ES06,S07}, S. Cuccagna, 
D. Pelinovsky and V. Vougalter \cite{CPV05}.
The eigenfunction expansions of $J$-selfadjoint operators were
not justified previously.

The spectral resolution
of bounded  $J$-selfadjoint nonnegative
operators in the Krein spaces
was constructed by M. Krein, H. Langer and Yu. Smul'jan \cite{KL1963,KS1966},
and extended to unbounded {\it definitizable} operators by M. Krein,
P. Jonas, H. Langer and others 
\cite{ IKL1982,  Jonas1981, Jonas1988,   L1981, LN1983}.
The corresponding unitary operators were considered by 
P. Jonas \cite{Jonas1986}.
However, the spectral resolution  alone is insufficient
for a justification of eigenfunction expansions. 
Our version 
of the theory under condition (\re{P})
allows us to justify the eigenfunction 
expansion (\re{AB0i}).

The spectral 
theory of 
definitizable operators was applied to the Klein-Gordon equations
with non-positive energy
by P. Jonas, H. Langer, B. Najman and C. Tretter
\cite{Jonas1993, Jonas2000, LNT2006,LNT2008,LT2006},
where the existence and uniqueness
of classical solutions were proved,
and 
the 
existence of unstable eigenvalues (imaginary frequencies) 
was studied.
The instability is related to the known {\it Klein paradox} 
in quantum mechanics \cite{Sakurai}.

The scattering theory 
for the  Klein-Gordon equations 
with non-positive energy 
was developed by
C. G\'erard and T. Kako 
using the theory of the 
definitizable operators in the Krein spaces 
\cite{Gerard2012, Kako1976}.

\medskip

The plan of our paper is as follows.
In Section 2 we 
justify the M. Krein transformation 
under condition (\re{P}),
and  construct a unitary dynamical group
and 
its  spectral representation.
In Section 3 we check condition (\re{P}) for the operator 
(\re{AB0i}), and 
in Sections 4-6
we justify
the eigenfunction expansion (\re{uns30i}) 
applying the methods of Section 2. 
In Appendix A we check the spectral condition (\re{P})
for a class of operators, 
and in Appendix B we check that $JB\ne BJ$ for 
linearizations  of nonlinear Schr\"odinger equations.
\medskip

{\bf Acknowledgments} The authors thank V. Ivrii, A. Kostenko,
M. Malamud 
and G. Teschl
for useful discussions on
pseudodifferential operators and $J$-selfadjoint operators.

\setcounter{equation}{0}
\section{Spectral representation}
We are going to obtain a spectral representation for solutions to 
equation (\re{NHl}).

\subsection{Generalized solutions}
Let  $D(B)$ denote the dense domain of the selfadjoint operator $B$.
We set $\Lambda:=B^{1/2}> 0$ and denote by
$\cV\subset \cX$
the Hilbert space
which is the domain  of  $\Lambda$ endowed with the norm
\begin{equation}\la{Vn22}
\Vert X\Vert_{\cV}:=\Vert \Lambda X\Vert_\cX~
\end{equation}
which is positive definite by (\re{P}).
We have a continuous injection of Hilbert spaces
$\cV\subset \cX$, and
the unitary operator
\begin{equation}\la{LL}
\Lambda: \cV\to \cX.
\end{equation}
In particular,
\begin{equation}\la{LR}
\Lambda^{-1}: \cX \to \cV~
\end{equation}
is a bounded operator.
For example, $\cV$ is the Sobolev space $H^1(\R^n)$
in the case of $\cX=L^2(\R^n)$ and 
operator $A= i(-\Delta+m^2)$  with any
real $m\ne 0$.

Finally, $J$ is the unitary operator in $\cX$ by (\re{Hs}),
and hence
\begin{equation}\la{AJV}
A^{-1}=-B^{-1}J: \cX \to \cV~
\end{equation}
is a bounded operator by (\re{LR}).
We will consider solutions
\begin{equation}\la{YD}
X(t)\in C(\R,\cV)
\end{equation}
to equation (\re{NHl}).
We will understand the equation in the sense of {\it mild solutions}
\cite{Cazenave}
\begin{equation}\la{ini}
X(t)-X(0)=A\int_0^t X(s)\5 ds~,~~~~~~~~~~~~t\in\R~.
\end{equation}
By (\re{AJV})
this is equivalent to
 the identity
\begin{equation}\la{ini2}
A^{-1}[X(t)-X(0)]= \int_0^t X(s)\5 ds~,
~~~~~~~t\in\R~,
\end{equation}
where the Riemann integral converges in $\cV$ by (\re{YD}).

\subsection{Krein transformation}

We apply well known  formal similarity transformation 
\begin{equation}\la{tra}
A=J\Lambda^2\to \Lambda A\Lambda^{-1}=\Lambda J\Lambda~, 
\end{equation}
rising to M. Krein.
This transformation corresponds to  the substitution
\begin{equation}\la{ZY}
Z(t):=\Lambda X(t)\in C(\R,\cX)~,
\end{equation}
used by M. Krein in 
the theory of parametric resonance:
see formula  (1.40) of
\cite[Chapter VI]{GK}.
Applying the  transformation (\re{tra}) to 
the both sides of
(\re{NHl}), we obtain formally the corresponding 'Schr\"odinger equation'
\begin{equation}\la{CPF3}
 i\dot Z(t)=H
Z(t)~,\,\,\,~~~~t\in\R~.
\end{equation}
Here $H$ stands for the `Schr\"odinger operator' 
\begin{equation}\la{CH}
H=\Lambda iJ\Lambda~, 
\end{equation}
which is symmetric  on the domain
\begin{equation}\la{DDH}
D(H)=\{X\in\cV\!:~ J\Lambda X\in\cV\}~.
\end{equation}
These arguments give all solutions to (\re{NHl}) in the case of finite dimensional
space $\cX$. In the  infinite dimensional case
the problem is less trivial.

\subsection{Selfadjoint generator}
We must justify that the operator $H$ is densely defined and selfadjoint
in our situation.
Otherwise, the formal transformation (\re{tra}) would not help 
to construct solutions of equation (\re{NHl}).

\bl\la{lGH}
Let condition \eqref{P} hold.
Then $D(H)$ is dense in $\cX$ and
the operator $H$ is
selfadjoint.
\el
\Pr
The operator $H$ is injective.
Further, $\Ran\Lambda=\cX$
by (\re{P}), and
 $J:\cX\to\cX$ is the unitary operator.
Hence,
$\Ran H=\cX$. Consider the inverse operator
\begin{equation}\la{G}
R:=H^{-1}=\Lambda^{-1} iJ\Lambda^{-1}~.
\end{equation}
It is selfadjoint since $D(R)=\Ran H=\cX$ and
$R$ is bounded and symmetric.
Finally, $R$ is injective, and hence, 
\begin{equation}\la{Hdef}
H=R^{-1}
\end{equation}
is a densely defined selfadjoint operator
by Theorem 13.11 (b) of \cite{Rudin}:
$$
H^*=H~, \quad D(H)=\Ran R\subset \cV~.
$$
Thus, the lemma is proved.
\bo
\medskip

As a result, 
under condition {\rm (\re{P})} the identity 
\begin{equation}\la{me2}
A=-i\Lambda^{-1}H\Lambda~
\end{equation}
holds on the domain $\Lambda^{-1}D(H)$, which is dense in $\cV$.
We will understand equation (\re{CPF3}) similarly to (\re{ini2}):
\begin{equation}\la{ini3}
iH^{-1}[Z(t)-Z(0)]= \int_0^tZ(s)\5 ds~,
~~~~~~~t\in\R~,
\end{equation}
where the Riemann integral converges in $\cX$ by (\re{ZY}).
\bc\la{ceqv}
Let condition {\rm (\re{P})} hold. Then
for any $Z(0)\in\cX$
equation {\rm (\re{CPF3})}
admits a unique solution $Z(t)\in C(\R,\cX)$ in the sense {\rm (\re{ini3})}.
The solution is given by
\begin{equation}\la{uns222}
Z(t)=e^{-iHt}Z(0)\in C(\R,\cX)~.
\end{equation}
\ec

\subsection{Unitary dynamical group}

It is easy to check that
equation
{\rm (\re{CPF3})} for $Z(t)\in C(\R,\cX)$
in the sense {\rm (\re{ini3})}
is equivalent to equation {\rm (\re{NHl})} in the sense (\re{ini2})
 for 
\begin{equation}\la{eq2}
X(t)=\Lambda^{-1}Z(t)\in C(\R,\cV)~.
\end{equation}
Hence, Corollary \re{ceqv} implies the following lemma.

\bl\la{cV} Let condition {\rm (\re{P})} hold. Then
\medskip\\
i) For any $X(0)\in\cV$, the function
\begin{equation}\la{uns3}
X(t)=\Lambda^{-1}e^{-iHt}\Lambda X(0)\in C(\R,\cV)~
\end{equation}
is the unique solution to {\rm (\re{NHl})} in the sense {\rm (\re{ini2})}.
\medskip\\
ii)
The dynamical group $V(t): X(0)\mapsto X(t)$ is unitary in $\cV$, since
\begin{equation}\la{UV}
\Vert X(t)\Vert_\cV:=\Vert \Lambda X(t)\Vert_\cX
=\Vert e^{-iHt} \Lambda X(0)\Vert_\cX
=\Vert \Lambda X(0)\Vert_\cX=:\Vert X(0)\Vert_\cV~.
\end{equation}
\el

\subsection{Spectral resolution}
Let $E(\omega)$ denote the spectral resolution of the selfadjoint operator $H$; i.e., 
\beqn
HZ=\int \omega\5 dE(\omega)Z~,~~Z\in D(H)~,\la{sreH2}
\eeqn
the Riemann--Stieltjes integrals converging in $\cX$.
This resolution and (\re{uns3}) imply
the following proposition.

\bp\la{tspec-c}
Let condition {\rm (\re{P})} hold. Then
the following spectral representation holds
for solutions to the complex Hamilton equations {\rm (\re{NHl})}:
\begin{equation}\la{srepc}
X(t)=\Lambda^{-1}\int e^{-i\omega t}dE(\omega)\Lambda X(0),~~X(0)\in\cV.
\end{equation}

\ep

Now let us apply our results to
real Hamilton equations (\re{NHlr}).
Let $B_r$ be a nonnegative invertible
selfadjoint operator on a~real Hilbert space $\cX_r$.
Then condition (\re{P}) holds for $B_r$.
\medskip\\
Let us consider complexifications $B$ and $J$ of the operators $B_r$ and $J_r$
 in the space $\cX=\cX_r+ i\cX_r$; these are defined by 
\begin{equation}\la{comp}
B(X_1+iX_2)=B_rX_1+iB_rX_2~,~~~~~J(X_1+iX_2)=J_rX_1+iJ_rX_2~,~~~~~~~~~~
X_1,~X_2\in\cX_r~.
\end{equation}
We claim that $B$ satisfies condition (\re{P}).
First, $B$ is invertible
since the complexification of $B_r^{-1}$ gives $B^{-1}$.
Moreover, $B^{-1}$ is a bounded
injective selfadjoint operator, and hence
$B$ is a
densely defined selfadjoint operator in $\cX$ by Theorem 13.11 (b) of~\cite{Rudin}.
Finally, $B$ is obviously nonnegative.

Further, let $\cV_r$ denote the domain of the real selfadjoint operator
$\Lambda_r:=B_r^{1/2}$.
Therefore,
 $\cV:=D(\Lambda)=\cV_r+i\cV_r$~, 
since $\Lambda:=B^{1/2}$ is the complexification of $\Lambda_r:=B_r^{1/2}$.
Hence, Corollary \re{cV} implies that equation
 (\re{NHl}) admits a unique solution
\begin{equation}\la{usol}
X(t)=Y_1(t)+iY_2(t)\in C(\R,\cV)~,~~~~~~~~~~Y_1(t), Y_2(t)\in \cV_r
\end{equation}
for any $X(0)=X_r\in \cV_r$~. 
The solution admits a spectral representation of type (\re{srepc}), where
$E(\omega)$ is the spectral family of the corresponding
 selfadjoint operator $H$ defined by~(\re{G}).
  Summarizing,
 Corollary \re{cV} and
Proposition \re{tspec-c} imply the following corollary. 
Let us denote by
$E_{rr}(\omega)$ and $E_{ir}(\omega)$ the real and imaginary
components of $E(\omega)|_{\cX_r}$~.

\bc\la{tspec-r}
Let $B_r$ be a nonnegative invertible
selfadjoint operator in a real Hilbert space $\cX_r$.
Then
\medskip\\
i) A solution $X_r(t)\in C(\R,\cV_r)$
to {\rm (\re{NHlr})} exists and is unique for any initial state $X_r(0)\in \cV_r$~,
and the energy {\rm (\re{hamr})} is conserved.
\medskip\\
ii) This solution
admits the
spectral representation  
\begin{equation}\la{srepr}
X_r(t)=\Lambda_r^{-1}\int
[\cos{\omega t}dE_{rr}(\omega)+\sin{\omega t}dE_{ir}(\omega)]\Lambda_r X_r(0)~,
~~X_r(0)\in\cV_r~,
\end{equation}
where
the Riemann--Stieltjes integral converges in $\cX_r$~.

\ec

\Pr
{\it i)} First, the real part $Y_1(t)\in C(\R,\cV_r)$ of solution (\re{usol})
is the unique solution to (\re{NHlr}) with $Y_1(0)=X_r$, 
since the real subspace $\cX_r$ is invariant with respect to $A_r=J_rB_r$.
Similarly, the imaginary part $Y_2(t)\in C(\R,\cV_r)$ vanishes, since $Y_2(0)=0$,
and hence
the energy conservation for the real solution follows from (\re{UV}) and (\re{hamr}).
\medskip\\
{\it ii)}
Formula (\re{srepr}) 
 follows from
(\re{srepc}).
The convergence of the Riemann-Stieltjes integral of
(\re{srepr}) in $\cX_r$ follows from
the convergence in $\cX$ of the corresponding integral of
(\re{srepc}) since
the convergence of a sequence in the Hilbert space $\cX=\cX_r+ i\cX_r$ is equivalent
to the convergence of the corresponding
 real and imaginary parts in $\cX_r$.
\bo

\setcounter{equation}{0}
\section{Application to  eigenfunction expansion}
We are going to apply 
our results
 to justification of the eigenfunction
expansion (\re{uns30i})
in the context of the system considered in \cite{KK11}.
We have used this expansion for the calculation of the
Fermi Golden Rule \cite[(5.14)]{KK11}.

\subsection{Linearization at the kink}
In \cite{KK11,KKm11} we studied the
 1D relativistic Ginzburg--\allowbreak Landau equation
\begin{equation}\la{GL}
\ddot\psi(x,t)=\fr{d^2}{dx^2}\psi(x,t)+F(\psi(x,t))~,\qquad x\in\R
\end{equation}
for real solutions $\psi(x,t)$.
Here 
$F(\psi)=-U'(\psi)$, where
$U(\psi)$ is similar to the
Ginzburg--\allowbreak Landau potential $U_{GL}(\psi) = (\psi^2-1)^2/4$, 
which corresponds to the cubic equation with $F(\psi) = \psi- \psi^3 $.
Namely, $U(\psi)$ is a real smooth even function satisfying the following
conditions:
\begin{equation}\la{U1}
U(\psi) >0~,~~~ \psi \ne \pm a~;~~~~~~~~~~~~~~~
U(\psi) = \fr{m^2}2(\psi\mp a)^2 + O(|\psi\mp a|^{14})~,~~~~ x\to \pm a~,
\end{equation}
where $a,m>0$.
The main goal of \cite{KK11,KKm11} was to prove of the asymptotic stability
of solitons (kinks) $\psi(x,t)=s_v(x-vt)$ that 
move with constant velocity $|v|< 1$, and
\begin{equation}\la{kink}
s_v(x)\to\pm a~,\qquad x\to\pm\infty~.
\end{equation}
Substituting  $\psi(x,t)=s_v(x-vt)$ into (\re{GL}), we obtain
\be\la{GLs}
v^2 s_v''(x)=s_v''(x)+F(s_v(x))~,\qquad x\in\R~.
\ee
The linearization
of (\re{GL})
at the kink $s_v(x-vt)$ in the moving frame
reads as (\re{NHl}) with $X=(\psi,\dot\psi)\in L^2(\R)\otimes\C^2$
(for the corresponding complexification)
and with the generator
\cite[(4.6)]{KKm11}
\begin{equation}\la{AA1}
A_v=\left(
\ba{cc}
v\fr d{dx}   & 1 \\
\fr{d^2}{dx^2}-m^2-V_v(x) & v\fr d{dx}
\ea
\right).
\end{equation}
Here the potential
\begin{equation}\la{Vv}
V_v(x)=-F'(s_v(x))-m^2\in C^\infty(\R)~.
\end{equation}
The kink $s_v(x)$ is an odd monotonic function
in a suitable coordinate $x$,
while $F'(\psi)=-U''(\psi)$ is an even function of $\psi$.
Hence, the potential $V_v(x)$ is the even function of $x$.
Moreover,
\begin{equation}\la{Vvm}
|V_v(x)|\le Ce^{-\kappa|x|}~,\qquad x\in\R~,
\end{equation}
where $\kappa>0$.
The generator (\re{AA1})
has the form $A_v=JB_v$ with
\begin{equation}\la{BJ}
B_v=\left(
\ba{cc}
S_v & -v\fr d{dx} \\
  v\fr d{dx}   &   1
\ea
\right)~,
\quad\quad
S_v:=-\fr{d^2}{dx^2} +m^2+V_v(x)
~,
\quad\quad
J:=\left(
\ba{cc}
0 & 1 \\
-1 & 0
\ea
\right).
\end{equation}
Obviously, $JB_v\ne B_vJ$.
Differentiating (\re{GLs}), we obtain 
\be\la{grst}
\ti S_v s_v'(x)=0~,~~~~~~~\ti S_v:=-(1-v^2)\fr{d^2}{dx^2} +m^2+V_v(x)~.
\ee



\subsection{Spectral condition}

Here we check condition (\re{P}) for operator $B_v$
in the space of the odd states in the case $v=0$.
The general case $|v|<1$ is considered in 
Appendix A (see Corollary \re{cP}).
We will write $A$,  $B$ and $S$ instead of   $A_0$, $B_0$ and $S_0$:
\begin{equation}\la{AB0}
A=\left(
\ba{cc}
0   & 1 \\
-S & 0
\ea
\right)~,~~~~~~~
B=\left(
\ba{cc}
S & 0 \\
0 & 1
\ea
\right)~,~~~~~~~~~~~~S:=-\fr{d^2}{dx^2}
+m^2+ V_0(x)~.
\end{equation}
The operators $B$ and $S$ 
are essentially selfadjoint
in $L^2(\R)\otimes\C^2$ and
$L^2(\R)$ respectively,
 by (\re{Vvm})
and
Theorems X.7 and X.8 of \cite{RS2}.
Now  (\re{grst}) with $v=0$ means 
that
$\lam=0\in \si_{pp}(S)$.
Moreover, $\lam=0$
is the minimal eigenvalue
of $S$, since the corresponding eigenfunction $s_0'(x)$
does not vanish \cite[(1.9)]{KKm11}.
Hence,
\begin{equation}\la{sB}
\sigma(S)\subset [0,\infty)~.
\end{equation}
Moreover, 
\begin{equation}\la{sB2}
\Ker S=(s_0'(x))~,
\end{equation}
since any second 
linearly independent
solution of the homogeneous equation
cannot belong to
$L^2(\R)$ by Theorem X.8 of \cite{RS2}.

Below we
restrict the phase space $L^2(\R)\otimes\C^2$
to the invariant subspace of the odd states
\begin{equation}\la{odd}
\cX=\{\psi\in L^2(\R)\otimes\C^2~:~ \psi(-x)=-\psi(x)~,~~~ x\in\R\}~, 
\end{equation}
as in \cite{KK11}.
The subspace is invariant for the operators $A$ and $B$, 
since the potential $V_0(x)$
is an even function.
Respectively, we consider the operator $S$ on the Hilbert 
space of odd functions
$L^2_{\rm odd}(\R)$.


\bl\la{lSB}
Condition {\rm (\re{P})} holds
for the operator $B$ on $\cX$.
\el
\Pr 
The ground state $s_0'(x)$ 
of $S$ on $L^2(\R)\otimes\C^2$
is an even function. Hence,
$\Ker S=0$
for $S$ on $\cX$ by (\re{sB2}).
Further, 
the continuous spectrum of $S$ lies in $[m^2,\infty)$. Hence,
(\re{sB}) implies that
\begin{equation}\la{spB0}
\sigma(S)=
 \{\lambda_1~,~....~,\lambda_N\}\cup [m^2,\infty)~,
\end{equation}
where $0<\lambda_1<~...~<\lambda_N<m^2$.
Finally, $\sigma(B)=\sigma(S)\cup \{1\}$, which implies (\re{P}).\bo
\medskip

We will assume below
the following spectral condition 
on the edge point of the continuous spectrum of $S$
imposed in \cite{KK11}:
\begin{equation}\la{SC11}
\mbox{\it The point $m^2$ is neither an eigenvalue
nor a~resonance of $S$~.}
\end{equation}

\setcounter{equation}{0}
\section{Orthogonal eigenfunction expansion}
Let us apply 
Proposition \re{tspec-c}
to the case of operators (\re{AB0}).
We have
\begin{equation}\la{LH0}
\Lambda:=B^{1/2}=\left(
\ba{cc}
\sqrt{S} & 0 \\
  0   &   1
\ea
\right),~~~~~~~~~~H:=\Lambda iJ\Lambda=
i\left(
\ba{cc}
0 & \sqrt{S} \\
 -\sqrt{S}  & 0
\ea
\right)=iJ\sqrt{S}~.
\end{equation}
Thus, 
$H$ is obviously selfadjoint in accordance to
 Lemma \re{lGH}.
Hence, (\re{spB0}) implies that
\begin{equation}\la{spB}
\sigma(H)
=(-\infty,-m]\cup
 \{\omega_{-N}~,~....~,\omega_{-1}~, \omega_1~,~...~,\omega_N\}\cup [m~,\infty)~,\qquad 
\omega_{\pm k}^2=\lambda_k~,\quad k=1~,~...~,N~.
\end{equation}
Respectively, formula (\re{srepc})
for solutions to (\re{NHl}) reads
\begin{equation}\la{uns30}
X(t)=
\Lambda^{-1}
\int_{\sigma(H)} e^{-i\omega t} d E(\omega)\Lambda X(0)
=
\sum_{-N}^N e^{-i\omega_k t} C_k a_k+
\Lambda^{-1}
\int_{\sigma_c(H)} e^{-i\omega t} d E(\omega)\Lambda X(0)~,
\end{equation}
where $\sigma_c(H)=\sigma_c=(-\infty, -m]\cup[m,\infty)$ is the
continuous spectrum of $H$, and
$a_k=
\Lambda^{-1} h_k
\in\cX$, where
$h_k$
are the eigenfunctions
of $H$ corresponding to the eigenvalues
$\omega_k$~.
Formula (\re{me2}) implies that $a_k$ are the eigenfunctions of $A$
corresponding to the eigenvalues
$-i\omega_k$~.

Let us denote by $X^c(t)$ the last integral
in (\re{uns30}):
\be\la{eife0}
X^c(t)
=
\Lam^{-1}\int_{\sigma_c}  e^{-i\omega t} d E(\omega)\Lambda X(0)~.
\ee
To prove (\re{uns30i})
it remains to justify the eigenfunction expansion 
\begin{equation}\la{eife}
X^c(t)
=
\int_{\sigma_c} e^{-i\omega t} C(\omega)\5 a_\omega\5 \5 d\omega~, 
\end{equation}
where $a_\omega$ are the
 generalized odd eigenfunctions from
the continuous spectral space of $A$
corresponding to the eigenvalues $-i\omega$~.
Then 
(\re{uns30i}) will follow from (\re{uns30}).

By  (\re{uns30}),  $X^c(t)$
is the solution to (\re{NHl}), and hence $Z^c(t):=\Lambda X^c(t)$ is the solution
to (\re{CPF3}).
We will deduce
(\re{eife}) from
 the corresponding representation
\begin{equation}\la{eife2}
Z^c(t)
=
\int_{\sigma_c} e^{-i\omega t} C(\omega)\5 h_\omega\5 \5 d\omega~, 
\end{equation}
where
$h_\omega$ are the generalized odd eigenfunctions
 of $H$ with the eigenvalues $\omega$.
We will prove (\re{eife2}) by solving equation (\re{CPF3})
for $Z^c(t)=(Z_1^c(t),Z_2^c(t))$.
By (\re{LH0}),
this equation is equivalent
to
the system
\begin{equation}\la{Z2}
\dot Z_1^c(t)=\sqrt{S}Z_2^c(t)~, \qquad 
\dot Z_2^c(t)=-\sqrt{S}Z_1^c(t)~.
\end{equation}
Eliminating $Z_2^c$, we obtain
\begin{equation}\la{Z3}
\ddot Z_1^c(t)= -SZ_1^c(t)~.
\end{equation}
Further we 
apply Theorem XI.41 of \cite{RS3} and
the arguments of \cite[pp 114-115]{RS3}. Namely,
the rapid decay (\re{Vvm}) and our spectral condition (\re{SC11}) imply
the following Limiting Absorption Principle (LAP) \cite{A,KKW12,RS3}:
\begin{equation}\la{LAP}
R(\lambda\pm i\varepsilon)\to R_\pm(\lambda),\qquad \varepsilon\to+0~,~~~~\lambda\in [m^2,\infty)~,
\end{equation}
where
$R(z):=(S-z)^{-1}$ and
the convergence holds
in the strong topology
of the space of continuous operators
$L^2_{s}\to L^2_{-s}$ with $s>1$.
Moreover, the traces of the resolvent $R_\pm(\lambda)$ are continuous functions of
$\lambda\ge m^2$
with values in $L(L^2_{s}, L^2_{-s})$, see \cite{A,KKW12}.
Here $L^2_{\rho}=L^2_{\rho}(\R)$ with $\rho\in\R$ denotes the
weighted Hilbert space with the norm
\begin{equation}\la{whs}
\Vert\psi\Vert^2_{L^2_{\rho}}:=\int \langle x \rangle^{2\rho}|\psi(x)|^2dx~,\qquad 
\langle x \rangle:=(1+x^2)^{1/2~}.
\end{equation}
The LAP serves as the basis for
the eigenfunction expansion 
\begin{equation}\la{eife3}
Z_1^c(t)=\int_{\sigma_c} d \cE(\omega^2)
[Z_1^c(0)\cos \omega t+Z_2^c(0)\sin \omega t]=
\int_{\sigma_c} e^{-i\omega t} C(\omega)\5 e_\omega\5 \5 d\omega~,
\end{equation}
where  $d\cE(\lambda)$ is the spectral resolution of $S$ in the space
$L^2_{\rm odd}(\R)$,
while
$e_\omega\in L^2_{-s}$
are generalized odd eigenfunctions
of $S$ corresponding to the eigenvalues $\omega^2\ge m^2$.
Here the first identity follows by Spectral Theorem, while the second
follows 
by  Theorem XI.41 (e) of \cite{RS3}.
The eigenfunctions are defined by formulas
of \cite[pp 114-115]{RS3}:
\begin{equation}\la{eigf}
e_\omega=W^*(\omega)f_\omega~,~~f_\omega(x):=\sin|\omega| x~,
~~ W(\omega):=[1+VR_0(\omega^2+i0)]^{-1}~,
~ \omega\in \sigma_c~.
\end{equation}
where $R_0(\lam):=(-\De+m^2-\lam)^{-1}$.
The operator $W(\omega)$
 is a~continuous function of $\omega\in\sigma_c$
with values in
 $L(L^2_{s},L^2_{s})$ by the formula 
\be\la{LL2}
[1+VR_0(\lam)]^{-1}=1-VR(\lam)
\ee
and the decay (\re{Vvm}).
Respectively,
the adjoint operator $W^*(\omega)$
is a continuous function of $\omega\in \sigma_c$
with values in
 $L(L^2_{-s},L^2_{-s})$.
As the result, $e_\omega$ is a
continuous function of $\omega\in\sigma_c$
with values in
$L^2_{-s}$. The normalization of $e_\omega$ coincides with the same of
the 'free' generalized eigenfunctions
$f_\omega$:
\begin{equation}\la{norm}
\langle e_{\omega}~, e_{\omega'}\rangle=\pi\5\delta(|\omega|-|\omega'|)~,
\quad \omega,\omega' \in\sigma_c
~,
\end{equation}
which follows from the last formula on page 115 of \cite{RS3}.
Finally,
Theorem XI.41 (e) of
\cite{RS3} implies that
the last integral (\re{eife3}) converges in
$L^2=L^2(\R)$:
\begin{equation}\la{icon0}
\Vert Z_1^c(t)-\int_{m\le |\omega|\le M} e^{-i\omega t} C(\omega)
\5 e_\omega\5 
\5 d\omega\Vert_
{L^2}
\to 0~,\qquad M\to\infty~.
\end{equation}

\br\la{rver}
Our modification
of the eigenfunction expansion theory
 differs  from {\rm \cite[pp 114-115]{RS3}} only in the use
of Hilbert spaces with weights
$\langle x \rangle^{s}$ instead of $e^{s |x|}$.
All conclusions of Theorem XI.41 from {\rm \cite{RS3}}
remain valid in the modified theory.
This modification allows us to apply Lemma \re{lambda}
below.

\er

Now (\re{eife2}) for $Z_1^c(t)$ follows from (\re{eife3}).
For $Z_2^c(t)$ we  use the first
equation of (\re{Z2}),  which implies
\beqn\la{eife4}
Z_2^c(t)
=
-i
\int_{\sigma_c} \sgn\omega \5\5\5e^{-i\omega t} C(\omega)\5 e_\omega\5 \5 d\omega~.
\eeqn
Combining (\re{eife3}) and (\re{eife4}), we obtain (\re{eife2})
with
\begin{equation}\la{Zk}
h_\omega:=\left(\ba{c} 1\\-i ~
\sgn\omega\ea\right)e_\omega \in L^2_{-s}\otimes \C^2~,\qquad 
\omega\in \sigma_c~,~~~~~~~\forall s>1~.
\end{equation}
Normalization (\re{norm}) implies the corresponding
normalization for $h_\omega$:
\begin{equation}\la{norm2}
\langle h_{\omega}, h_{\omega'}\rangle=2\pi\5\delta(\omega-\omega')~,
\qquad \omega,\omega' \in\sigma_c
~.
\end{equation}
\bl\la{lor} Let condition (\re{SC11}) hold and $s>1$. Then 
\medskip\\
i) $h_\omega$ are generalized eigenfunctions of $H$, i.e., 
\begin{equation}\la{eigf2}
H Z^c(t)=\int_{\sigma_c} e^{-i\omega t} \omega\5
C(\omega)\5 h_\omega\5 \5 d\omega~~~~~~for~~~~Z^c(t)\in D(H)~.
\end{equation}
\medskip\\
ii) $h_\omega$ is a continuous function of $\omega\in\sigma_c$ with values in
 $L^2_{-s}\otimes \C^2$~.
\medskip\\
iii)
The integral (\re{eife2})
converges in
$L^2\otimes\C^2$ in the following sense{\rm :}
\begin{equation}\la{icon}
\Vert Z^c(t)-\int_{m\le |\omega|\le M} e^{-i\omega t} C(\omega)\5 h_\omega
\5 \5 d\omega\Vert_{L^2\otimes\C^2}\to 0~, \qquad M\to\infty~.
\end{equation}

\el
\Pr i) $Z^c(t)\in D(H)$ means that $Z^c_{1,2}(t)\in D(\sqrt{S})$. 
Furthermore,
\begin{equation}\la{eigf3}
H Z^c(t)=i\sqrt{S}
\left(\ba{c}
 Z^c_2(t)\\ -Z^c_1(t)
\ea
\right).
\end{equation}
Now (\re{eigf2}) follows from the expansions (\re{eife3}) and (\re{eife4})
for $Z^c_{1,2}(t)$ by
\cite[Theorem XI.41 (c)]{RS3}, since $e_\omega$ are the generalized 
eigenfunctions of $S$
with the eigenvalues $\omega^2$, and {\it formally},
\begin{equation}\la{eigf4}
i\sqrt{S}
\left(\ba{c}
 -i~\sgn\omega\\ -1
\ea
\right)e_\omega=
\left(\ba{c}
\sgn\omega\\ -i
\ea
\right)|\omega|e_\omega=\omega h_\omega~.
\end{equation}
\medskip\\
ii) 
$h_\omega$ is a continuous function of $\omega\in\sigma_c$  by similar property of $e_\omega$~.
\medskip\\
iii) (\re{icon}) follows from (\re{icon0}) and similar convergence for
$Z_2^c$.\bo

\br
The generalized eigenfunctions 
(\re{Zk})
are 
proportional to  $e_\omega$
since $H$ and $S$ commute with each other.
This argument was the main idea of our derivation
of the  eigenfunctions  (\re{Zk}).

\er

\setcounter{equation}{0}
\section{Non-orthogonal eigenfunction expansion}
Let us denote by  $Z^c_M(t,x)$ the integral in (\re{icon}).
It  is defined 
for almost all $x$, i.e.,
\begin{equation}\la{icon4}
 Z^c_M(t,x):=
\int_{m\le |\omega|\le M} e^{-i\omega t} C(\omega)\5 h_\omega(x)
\5 \5 d\omega~, ~~~~~~~~~~~~a.a.~~~x\in\R~.
\end{equation}
To justify (\re{eife}) we should adjust the meaning of this integral
by the following lemma.

\bl\la{lsi}
Let condition {\rm (\re{SC11})} hold and $s >1$.
Then 
\medskip\\
i)
The integral {\rm (\re{icon4})} 
converges absolutely in $L^2_{-s}\otimes\C^2$:
\begin{equation}\la{whs2}
\int_{m\le |\omega|\le M} \Vert C(\omega)\5 h_\omega\Vert_{L^2_{-s}\otimes \C^2}\5 d\omega
<\infty~,
~~~~~~M>m~.
\end{equation}
ii) The integral  {\rm (\re{icon4})} coincides a.e.  with the corresponding integral of
the
$L^2_{-s}\otimes\C^2$-valued integrand.

\el
\Pr
i) 
 (\re{norm2}) and
the Plancherel identity \cite[(82e$'$)]{RS3} imply that
\begin{equation}\la{Pide}
\Vert\int_\alpha^\beta C(\omega)\5 h_\omega\5 \5 d\omega\Vert_{\cX}^2=
2\pi\int_\alpha^\beta |C(\omega)|^2\5 d\omega=
\Vert P_{[\alpha,\beta]}Z^c(0)\Vert^2_{\cX}\le \Vert Z(0)\Vert_\cX^2~,
~~~~~~[\alpha,\beta]\subset \sigma_c~.
\end{equation}
Hence, (\re{whs2}) follows
by the Cauchy--Schwarz inequality
and Lemma \re{lor} i).
\medskip\\
ii)
The scalar products of the both integrals with any test function $\varphi\in C_0^\infty(\R)$
coincide  
by the Fubini theorem  since $h_\omega(x)$ can be chosen a measurable function 
of $(\omega,x)\in \sigma_c\times\R$ by  Lemma \re{lor} ii).
\bo
\medskip

Now we are going to deduce (\re{eife}) applying $\Lambda^{-1}$ to
(\re{eife2}).

\bl\la{lambda}
The operator $\Lambda^{-1}:
L^2_{\rho}\otimes\C^2\to L^2_{\rho}\otimes\C^2$ is continuous for every
$\rho\in\R$.

\el
\Pr
By (\re{LH0}) we have
\begin{equation}\la{pi0}
\Lambda^{-1}=
\left(
\ba{cc}
S^{-1/2} & 0 \\
  0   &   1
\ea
\right)~.
\end{equation}
Hence,
it suffices to prove the continuity for
the operator
$S^{-1/2}$ in $L^2_{\rho}$, which means the
continuity of operator
$$
\langle x\rangle^{\rho} S^{-1/2}
\langle x\rangle^{-\rho}:~~ L^2(\R)\to L^2(\R)~.
$$
This continuity
follows by the Theorem of Composition of PDO, since $S^{-1/2}$
is a PDO 
of the class $HG^{-1,-1}_1$, see definition 25.2 in \cite{Shubin}.
This follows from
\cite[Theorem 29.1.9]{Hor4}
and also by
an extension of \cite[Theorem 11.2]{Shubin} to PDO
with nonempty continuous spectrum. It is important that
operator $S$ and its main symbol $\xi^2$ satisfy
$$
\xi^2\not\in (-\infty,0]~, ~~\xi\ne 0~;
~~~~~~\sigma(S)\cap (-\infty,0]=\emptyset~.
$$
Hence, conditions (10.1) and (10.2)
of \cite{Shubin} hold.
\bo
\medskip

Lemma \re{lambda} and (\re{Zk}) imply that
\begin{equation}\la{Xom}
a_\omega:=\Lambda^{-1}h_\omega\in L^2_{-s}\otimes\C^2~,~~~~~~~s>1~.
\end{equation}
\bl\la{leig} 
$a_\omega$ are generalized eigenfunctions of $A$ corresponding to
the eigenvalues
$-i\omega$.
\el
\Pr
Let $P_c^M$ denote the spectral projection $\chi_c^M(H)$, where
$\chi_c^M$ is the indicator of
$\sigma_c^M:=\sigma_c\cap[-M,M]$. 
Let us take any $Z^c(0)\in P^M_c\cX$.
Then $Z^c(0)\in D(H)$, and hence 
(\re{eife2}) and (\re{eigf2})
with $t=0$ imply
the expansions
\begin{equation}\la{fi}
Z^c(0)=\ds\int_{\sigma_c^M} C(\omega)\5 h_\omega\5 \5 d\omega~,
~~~~~~~~~~~~~~~~~~~
HZ^c(0)
=
\int_{\sigma_c^M} \omega \5\5 C(\omega)\5 h_\omega\5 \5 d\omega~.~~~~~~~~~~~~~~~~~~
\end{equation}
Now let us take any $X\in \Lambda^{-1}P_c^M\cX$ and
write (\re{fi})  for
 $Z^c(0)=\Lambda X$.
Applying $\Lambda^{-1}$ to each side, we obtain
\begin{equation}
\la{eife22r}
X=\ds\int_{\sigma_c^M} C(\omega)\5 a_\omega\5 \5 d\omega~,
~~~~~~~~
\Lambda^{-1} H\Lambda X=
AX=
-i\int_{\sigma_c^M} \omega \5\5 C(\omega)\5 a_\omega\5 \5 d\omega~,~~~~~~~~~~
\end{equation}
where we have used definition (\re{Xom}), Lemmas \re{lsi} and \re{lambda},~
and the expression (\re{me2}) for $A$.
Identities (\re{eife22r})
 mean that
$a_\omega$ are generalized eigenfunctions
 in the sense of \cite[(80b)]{RS3}.
\bo
\medskip

Finally, the main result of our paper is the following.

\bt\la{tmain}
Let condition {\rm (\re{SC11})} hold,
 $X(0)\in \cV$ and $s>1$. Then eigenfunction
expansion \eqref{eife} holds in the following sense{\rm :}
\begin{equation}\la{icon2}
\Vert X^c(t)
-\int_{m\le |\omega|\le M} e^{-i\omega t} C(\omega)\5 a_\omega\5 \5 d\omega\Vert_\cV\to 0~,
~~~~~~~~M\to\infty~,
\end{equation}
where the integral converges
in  $L^2_{-s}\otimes\C^2$, and hence, a.e., as in  {\rm (\re{icon4})}.

\et
\Pr
We should prove that
\begin{equation}\la{icon22}
\Vert\Lambda X^c(t)
-\Lambda\int_{m\le |\omega|\le M} e^{-i\omega t} C(\omega)\5 a_\omega\5 \5 d\omega\Vert_\cX\to 0~,
~~~~~~~~M\to\infty~.
\end{equation}
First, we recall that
$\Lambda X^c(t)=Z^c(t)$. 
Second,
\begin{equation}\la{icon3}
\int_{m\le |\omega|\le M} e^{-i\omega t} C(\omega)\5 a_\omega\5 \5 d\omega
=\Lambda^{-1}\int_{m\le |\omega|\le M} e^{-i\omega t} C(\omega)\5 h_\omega\5 \5 d\omega
\end{equation}
by definition (\re{Xom}) and Lemmas \re{lsi}, \re{lambda}.
Now (\re{icon22}) follows from (\re{icon}).
\bo

\setcounter{equation}{0}
\section{Symplectic normalization}
Now let us renormalize $h_\omega$ as follows:
\be\la{nor}
\langle h_{\omega}, h_{\omega'}\rangle=|\omega|\5\delta(\omega-\omega')~,~~~~~~~\omega,\omega' \in \sigma_c
~.
\ee
This means that for any $M<\infty$
\be\la{nor2}
\langle Z_1,Z_2 \rangle=\int_{m\le|\om|\le M}|\om|\5 C_1(\omega)\ov{C_2(\om)}d\om
~,~~~~~~~~\mbox{\rm for}~~~~~~Z_{1,2}=\int_{m\le|\om|\le M}C_{1,2}(\omega) h_\om d\om\in \cX~.
\ee
Let us express these formulas
in terms of
$X_{1,2}:=\Lambda^{-1}Z_{1,2}\in\cV$
and
the eigenfunctions $a_\omega:=\Lambda^{-1} h_\omega$. 

First, 
(\re{eigf2}) and  (\re{nor2})  imply that
\be\la{nor3}
\langle H^{-1}Z_1,Z_2 \rangle=
\int_{m\le|\om|\le M}\sgn\om\5\5 C_1(\omega)\ov{C_2(\om)}d\om~.
\ee
On the other hand, (\re{G}) implies that
\be\la{nor4}
\langle H^{-1}Z_1,Z_2 \rangle=\langle \Lam^{-1}iJ\Lam^{-1}Z_1,Z_2 \rangle
=-i\langle \Lam^{-1}Z_1,J\Lam^{-1}Z_2 \rangle
=-i\langle X_1,JX_2 \rangle~.
\ee
Finally, (\re{nor3}) - (\re{nor4}) imply that
\be\la{nor5}
-i\langle X_1,JX_2 \rangle=
\int_{m\le|\om|\le M}\sgn\om\5\5 C_1(\omega)\ov{C_2(\om)}d\om
~~~~~\mbox{\rm for}~~~~~
X_{1,2}=\int_{m\le|\om|\le M}C_{1,2}(\omega) a_\om d\om~.
\ee
By definition, this means that
\be\la{nor6}
\langle a_{\omega}~,~ Ja_{\omega'}\rangle=i\sgn\omega\5\5\5 \delta(\omega-\omega')~,
~~~~~~~~~~
\omega,\omega'\in \sigma_c~.
\ee
Now the expansion \eqref{eife}  coincides with \cite[(2.1.13)]{BS03}, thereby justifying our calculation
of the Fermi Golden Rule \cite[(5.14)]{KK11}.



\appendix

 \setcounter{section}{0}
\setcounter{equation}{0}
\protect\renewcommand{\thesection}{\Alph{section}}
\protect\renewcommand{\theequation}{\thesection. \arabic{equation}}
\protect\renewcommand{\thesubsection}{\thesection. \arabic{subsection}}
\protect\renewcommand{\thetheorem}{\Alph{section}.\arabic{theorem}}

\section{Spectral condition for $|v|<1$}
Let us check the spectral condition (\re{P}) for operators
$B_v$ from (\re{BJ}) with any $|v|<1$ in the space of the odd states.
First let us check the continuous spectrum.
\begin{lemma}\la{lGL}
The continuous spectrum of  $B_v$ lies in 
$[\de,\infty)$
with some $\de>0$.
\end{lemma}
\Pr
By 
Corollary 2 (c) of \cite[XIII.4]{RS4}
and (\re{Vvm})
it suffices to find the continuous spectrum of the unperturbed
operator $B_v^0$ corresponding to $V_v(x)=0$.
Consider the spectral equation
\be\la{mice}
 (B_v^0-\lam)\psi=0
\ee
and find the solution of type $\psi=e^{ikx}\phi$ with real $k$ and $\phi\in\C^2$.
Substituting to (\re{mice}) we obtain
$$ 
\left(
\ba{cc}
k^2+m^2-\lam  & ikv \\
-ikv & 1-\lam
\ea
\right)\phi=0~.
$$
For nonzero vectors $\phi$, the determinant of the matrix vanishes:
\be\la{det}
  k^2(1-\lam-v^2)+(m^2-\lam)(1-\lam)=0~.
\ee
Then $k^2=(m^2-\lam)(1-\lam)/(\lam-1+v^2)\ge 0$. This inequality holds if
$$
\left\{
\ba{lll} 
\lam\in [1-v^2,1]\cup[m^2,\infty) &{\rm for}&\quad 1\le m^2~,\\ \\
\lam\in [m^2,1-v^2)\cup (1,\infty) &{\rm for}&\quad m^2\le 1-v^2\le 1~,\\ \\
\lam\in [1-v^2,m^2]\cup[1,\infty) &{\rm for}&\quad 1-v^2\le m^2\le 1~.
\ea
\right.
$$
Now the lemma follows with  $\de =\min(1-v^2,m^2)$.
\bo
\medskip

Now let us consider the discrete spectrum.

\begin{lemma}\la{GL1}
i) $\dim\Ker B$ is generated by $(s_0'(x),-vs_0''(x))$, where $s_0'(x)$ is an even function.
\\
ii)
 The nonzero discrete spectrum of $B$ is positive.
\end{lemma}
\Pr
{\it i)}
Equation  $B\psi=0$ is equivalent to the system
\be\la{dsBv0}
\left(
\ba{cc}
S_v &   -v\fr d{dx} \\
     v\fr d{dx} \         &         1
\ea
\right)\left(
\ba{c}
\psi_1 \\ \psi_2
\ea
\right)=0~.
\ee
The second equation (\re{dsBv}) implies
$\psi_2=-v \psi_1'$.
Substituting into the first equation we obtain
\be\la{Hv10}
\ti S_v \psi_1=0~,
\ee
where $\ti S_v$ is defined in (\re{grst}).
However, the equation (\re{grst}) with $v=0$ means 
that
$\lam=0\in \si_{pp}(\ti S_v)$.
Moreover, $\lam=0$
is the minimal eigenvalue
of $\ti S_v$, since the corresponding eigenfunction $s_v'(x)$
does not vanish \cite[(1.9)]{KKm11}.
Hence,
\begin{equation}\la{sBa}
\sigma(\ti S_v)\subset [0,\infty)~.
\end{equation}
Moreover, 
\begin{equation}\la{sB2a}
\Ker \ti S_v=(s_0'(x))~,
\end{equation}
since any second 
linearly independent
solution of the homogeneous equation
cannot belong to
$L^2(\R)$ by Theorem X.8 of \cite{RS2}.
\medskip\\
{\it Step ii)}
Consider equation $B_v\psi=\lam \psi$ with $\lam<0$:
\be\la{dsBv}
\left(
\ba{cc}
S_v-\lam &   -v\fr d{dx} \\
     v\fr d{dx} \         &         1-\lam
\ea
\right)\left(
\ba{c}
\psi_1 \\ \psi_2
\ea
\right)=0~.
\ee
The second equation (\re{dsBv}) implies
$\psi_2=v \psi_1'/ (\lam-1)~.$
Substituting into the first equation we obtain
\be\la{Hv1}
(\ti S_v+\fr{v^2\lam}{1-\lam}\fr{d^2}{dx^2}-\lam)\psi_1=0~.
\ee
For $\lam<0$ the operator is positive
since 
$\ti S_v\ge 0$ by (\re{sBa}).
Hence, equation (\re{Hv1}) has no nonzero solutions $\psi_1\in L^2$.
\bo
\bc\la{cP}
Lemmas \re{GL} and \re{GL1} imply that
 condition (\re{P})  holds for  $|v|<1$ in the space of the odd states.
\ec

\section{Linearization of $U(1)$-invariant Hamilton PDEs}

Equations (\re{NHl})
with $JB\ne BJ$
arise in the linearization of nonlinear
$U(1)$-invariant Hamilton PDEs.
Namely, consider the $U(1)$-invariant
Hamilton functional
\begin{equation}\la{NSH}
\cH(\psi)=\fr12\int\Big[|\nabla\psi(x)|^2+U(x,|\psi(x)|^2)\Big]dx
\end{equation}
with a real potential $U(x,r)$ and $\psi(x)\in\C=\R^2$.
The corresponding Hamilton equation reads as the nonlinear
Schr\"odinger equation
\begin{equation}\la{S}
i\dot \psi(x,t)
=\nabla_{\ov\psi}\cH(\psi)
=-\Delta\psi(x,t)+U_r(x,|\psi|^2)\psi~,~~~~~x\in\R^n~,
\end{equation}
where $i$ can be regarded as a~real $2\times 2$ matrix $J$
of type (\re{BJ}).
The linearization at a stationary state $s_0(x)$
is obtained by substitution $\psi=s_0+\varphi$ and expansion
$|\psi|^2=|s_0|^2+2s_0\cdot\varphi+|\varphi|^2$,
where $s_0\cdot\varphi$ is the scalar product of the real vectors
from $\R^2$. Neglecting the terms of higher order, we obtain
the linearized equation
\begin{equation}\la{Sl}
i\dot \varphi(x,t)=-\Delta\varphi(x,t)+ U_r(x,|s_0(x)|^2)\varphi+
2U_{rr}(x,|s_0(x)|^2)(s_0\cdot\varphi) s_0~, 
\end{equation}
which can be represented in the form (\re{NHl}) with $J=-i$
and $X(t)=(\rRe\vp(t),\rIm\vp(t))$.
The last term of (\re{Sl})
is not complex linear operator of $\varphi$.
In other words,
it
does not commute
with the multiplication of $\varphi$ by $i$. So $JB\ne BJ$
if $U_{rr}(x,|s_0(x)|^2)\not\equiv 0$.
Let us assume that $U(x,r)$ is a 
real-analytic function of $r>0$, and $s_0(x)\not\equiv 0$.
Then
the last term
of (\re{Sl})
vanishes exactly for the
 linear Schr\"odinger equation
when
$U(x,r)=V(x)r$.



\begin{thebibliography}{99}


\bibitem{A}
S. Agmon, Spectral properties of Schr\"odinger operator
and scattering theory,
{\em Ann. Scuola Norm. Sup. Pisa}, Ser.~IV {\bf 2}, 151-218 (1975).




\bibitem{AI1989}
T.Ya. Azizov, I.S. Iokhvidov,             
  Linear Operators in Space with an Indefinite Metric,
John Wiley \& Sons, Chichester, 1989.


      

\bibitem{BC11}
D. Bambusi, S. Cuccagna,
On dispersion of small energy solutions of the nonlinear Klein--Gordon equation with a potential,
{\em Amer. J. Math.} {\bf 133} (2011), no. 5, 1421-1468.

\bibitem{BC12}
N. Boussaid, S. Cuccagna,
On stability of standing waves of nonlinear Dirac equations,
{\em Comm. PDE} {\bf 37} (2012), no. 6, 1001-1056.
arXiv:1103.4452.



\bibitem{BP93} V.S. Buslaev, G.S. Perelman,
Scattering for the nonlinear
Schr\"odinger equation: states close to a soliton,
{\em St.Petersburg Math. J.} {\bf 4} (1993), 1111-1142.

\bibitem{BP95} V.S. Buslaev, G.S. Perelman,
On the stability of solitary waves
for nonlinear Schr\"odinger equations,
{\em Nonlinear evolution equations},
Transl. Ser. 2, 164, Amer. Math. Soc., Providence, RI, 1995, pp. 75-98.

\bibitem{BS03}
V.S. Buslaev, C. Sulem,
On asymptotic stability of solitary waves for nonlinear
Schr\"odinger equations,
{\em Ann. Inst. Henri Poincar\'e, Anal.
Non Lin\'eaire}
{\bf 20} (2003), no. 3, 419-475.
\bibitem{Cazenave}
T. Cazenave, A. Haraux,
Semilinear evolution equations, Clarendon Press, Oxford, 1998.


\bibitem{C01} S. Cuccagna,
Stabilization of solutions to nonlinear
Schr\"odinger equations,
{\em Comm. Pure Appl. Math.} {\bf 54} (2001), 1110-1145.

\bibitem{C03} S. Cuccagna,
On asymptotic stability of ground states of NLS,
{\em Rev. Math. Phys.} {\bf 15} (2003), 877-903.

\bibitem{CPV05}
S. Cuccagna, D. Pelinovsky,V. Vougalter,            
  Spectra of positive and negative energies in the linearized NLS problem,
{\em  Commun. Pure Appl. Math.} {\bf 58} (2005), no. 1, 1-29.
                              


\bibitem{C11}
S. Cuccagna,
On scattering of small energy solutions of non-autonomous
Hamiltonian nonlinear Schr\"odinger equations,
{\em J. Differ. Equations} {\bf 250} (2011), no. 5, 2347-2371.
\bibitem{C11H}
S. Cuccagna,
The Hamiltonian structure of the nonlinear Schr\"odinger equation and the asymptotic
stability of its ground states,
{\em Commun. Math. Phys.} {\bf 305} (2011), no. 2, 279-331.
\bibitem{C12}
S. Cuccagna,
On asymptotic stability of moving ground states of the nonlinear Schr\"odinger equation,
To appear in {\em Trans. Amer. Math. Soc}, 2012. arXiv:1107.4954




\bibitem{CM08}
S. Cuccagna, T. Mizumachi,
On asymptotic stability in energy space of ground states for nonlinear 
Schr\"odinger equations,
{\em Commun. Math. Phys.} {\bf 284} (2008), no. 1, 51-77.


\bibitem{CT09}
S. Cuccagna, M. Tarulli,
On asymptotic stability in energy space of ground states of NLS in 2D,
 {\em Ann. Inst. Henri Poincar\'e, Anal. Non Lin\'eaire}
 {\bf 26} (2009), no. 4, 1361-1386.


\bibitem{ES06}
M.B. Erdogan, W. Schlag, 
Dispersive estimates for Schrödinger operators in the presence of 
a resonance and/or an eigenvalue at zero energy in dimension three. II,
{\em J. Anal. Math.} {\bf 99}(2006), 199-248.
                              




\bibitem{Gerard2012}
C. G\'erard, 
 Scattering theory for Klein-Gordon equations with non-positive energy,
{\em Ann. Henri Poincar\'e} {\bf 13} (2012), no. 4, 883-941.



\bibitem{GK}
I.C. Gohberg, M.G. Krein,
Theory and applications of Volterra operators in Hilbert space,
AMS, Providence, R.I., 1970.

\bibitem{Hor4}
L. H\"ormander,           
The analysis of linear partial differential operators. IV: Fourier integral operators,
Springer,  Berlin, 2009.
  
\bibitem{IKL1982}
I.S. Iohvidov, M.G. Krein, H. Langer,            
  Introduction to the spectral theory of operators 
in spaces with an indefinite metric,
Mathematical Research, Vol. 9,  Akademie-Verlag, Berlin, 1982.
      



\bibitem{Jonas1981}
P. Jonas, 
 On the functional calculus and the spectral function for definitizable operators in Krein space
 {\em Beitr. Anal.} {\bf 16}(1981),  121-135.
                              

\bibitem{Jonas1986} 
P. Jonas,  On a class of unitary operators in Krein space, pp 151-172 in: 
Operator
Theory: Advances and Applications, Vol.17, Birkh\"auser Verlag, Basel-Boston-Stuttgart, 1986.

\bibitem{Jonas1988} 
P. Jonas,  On a class of selfadjoint operators in Krein space and their compact
perturbations,
{\em Integral Equ. Oper. Theory} {\bf 11} (1988), 351-384.


\bibitem{Jonas1993}
P.  Jonas,  On the spectral theory of operators associated with perturbed Klein-
Gordon and wave type equations,
{\em  J. Oper. Theory} {\bf 29} (1993), 207-224.

\bibitem{Jonas2000}
P.  Jonas,  On bounded perturbations of operators of Klein-Gordon type,
{\em  Glasnik
Math.} {\bf 35} (2000), 59-74.





\bibitem{Kako1976} Kako, T.: Spectral and scattering theory for the J-selfadjoint operators associated
with the perturbed Klein- Gordon type equations. J. Fac. Sci. Univ. Tokyo
Sec. I A 23, 199â221 (1976)


\bibitem{IKS12} V. Imaykin, A.I. Komech, H. Spohn,
Scattering asymptotics for a charged particle coupled to the Maxwell field,
{\em J. Math. Physics} {\bf 52} (2011), no. 4, 042701-042701-33.
arXiv:0807.1972


\bibitem{KKS11}
 A.I. Komech, E.A. Kopylova, H. Spohn,
Scattering of solitons for Dirac equation coupled to a particle,
{\em J. Math. Analysis and Appl.} {\bf 383} (2011), no.~2, 265--290. arXiv: 1012.3109



\bibitem{KK11}
 E.A. Kopylova, A.I. Komech,
On asymptotic stability of kink for relativistic Ginzburg--\allowbreak Landau
equation, {\em Arch. Rat. Mech. Anal.} {\bf 202} (2011),
no. 2, 213--245. arXiv:0910.5539

\bibitem{KKm11}
 E.A. Kopylova, A.I. Komech,
On asymptotic stability of moving kink for relativistic
Ginzburg--\allowbreak Landau equation,
{\em Comm. Math. Physics} {\bf 302} (2011), no.1,
225-252. arXiv:0910.5538


\bibitem{KKW12}
A. Komech, E.A. Kopylova,
Dispersion decay and scattering theory, Wiley, Hoboken, NJ, 2012.





\bibitem{KL1963}   M.G. Krein, H.K. Langer,          
  The spectral function of a selfadjoint operator in 
a space with indefinite metric,
{\em Sov. Math. Dokl.} {\bf 4} (1963), 1236-1239. 

\bibitem{KS1966}
M.G. Krein, Yu. Shmul'jan,
$J$-polar representations of plus-operators, {\em Mat. Issled.} {\bf 1} (1966), no.2, 172-210.
[Russian]



\bibitem{L1981}
H. Langer,
Spectral functions of definitizable operators in Krein spaces, 
pp. 1-46 in:
D. Butkovic, H. Kraljevic, S. Kurepa,
 Functional Analysis, LNM0948, Berlin, Springer, 1981.

\bibitem{LN1983}
Langer, H.; Najman, B.: Perturbation theory for definitizable operators
in Krein spaces, J.Operator Theory 9 (1983), 297-317.

\bibitem{LNT2006}
 H. Langer, B. Najman, C.  Tretter,          
  Spectral theory of the Klein-Gordon equation in Krein spaces,
{\em Proc. Edinb. Math. Soc., II. Ser.} 
{\bf 51} (2008), no. 3, 711-750.
                              

\bibitem{LNT2008}
H. Langer, B. Najman, C.  Tretter,  
Spectral theory of the Klein-Gordon equation in Pontryagin spaces,
{\em Commun. Math. Phys.} {\bf 267} (2006), no. 1, 159-180.
                             

\bibitem{LT2006}
 H. Langer, C.  Tretter,  
Variational principles for eigenvalues of the Klein-Gordon equation,
{\em J. Math. Phys.} {\bf 47} (2006), no. 10, 103506, 18 p.
                              



\bibitem{MM08} Y.~Martel, F.~Merle,
Asymptotic stability of solitons of the gKdV
equations with general nonlinearity,
{\em Math. Ann.} {\bf 341} (2008), 391-427.

\bibitem{MW96} J. Miller, M. Weinstein,
Asymptotic stability of solitary
waves for the regularized long-wave equation,
{\em Comm. Pure Appl. Math.} {\bf 49} (1996), 399-441.



\bibitem{PW94} R.L. Pego, M.I. Weinstein,
Asymptotic stability of solitary waves,
{\em Comm. Math. Phys.} {\bf 164} (1994), 305-349.


\bibitem{RS1} M.~Reed, B.~Simon,
 Methods of modern mathematical physics I:
 Functional Analysis, Academic Press, NY, 1980.

\bibitem{RS2} M.~Reed, B.~Simon,
 Methods of modern mathematical physics II:
 Fourier Analysis, Self-Adjointness, Academic Press, NY, 1975.


\bibitem{RS3} M.~Reed, B.~Simon,
 Methods of modern mathematical physics III:
 Scattering Theory, Academic Press, NY, 1979.

\bibitem{RS4} M.~Reed, B.~Simon,
 Methods of modern mathematical physics IV:
 Analysis of Operators, Academic Press, NY, 1978.



\bibitem{Rudin}
W. Rudin,
 Functional analysis, McGraw-Hill, New York, 1991.

\bibitem{Sakurai}
J.J.Sakurai, Advanced Quantum Mechanics, 
Addison-Wesley, Reading, Mass.,  1967.


\bibitem{S07}
W.Schlag,            
  Dispersive estimates for Schrödinger operators: a survey,
pp 255-285 in:
J. Bourgain (ed.) et al., 
Mathematical aspects of nonlinear dispersive equations. 
Lectures of the CMI/IAS workshop on mathematical aspects of nonlinear PDEs, 
Princeton, NJ, USA, 2004. 
 NJ: Princeton University Press, Princeton, 2007.


\bibitem{Shubin}
M.A. Shubin,
Pseudodifferential operators and spectral theory,
Springer, NY, 1987.



\bibitem{Seeley} R.T.~Seeley, 
Complex powers of an elliptic operator,
{\em Proc. Sympos. Pure Math.} {\bf 10} (1967), 288-307.

\bibitem{Sigal93} I.M. Sigal,
Nonlinear wave and Schr\"odinger equations.
I: Instability of periodic and quasiperiodic solutions,
{\em Commun. Math. Phys.} {\bf 153} (1993), no.2, 297-320.



\bibitem{SW90} A. Soffer, M.I. Weinstein,
Multichannel nonlinear scattering in
nonintegrable systems,
{\em Comm. Math. Phys.} {\bf 133} (1990), 119-146.

\bibitem{SW92} A. Soffer, M.I. Weinstein,
Multichannel nonlinear scattering and stability II.
The case of anisotropic and potential and data,
{\em J. Differential Equations} {\bf 98} (1992), 376-390.

\bibitem{SW99} A. Soffer, M.I. Weinstein,
Resonances, radiation damping and instability in
Hamiltonian nonlinear wave equations,
{\em Invent. Math.} {\bf 136} (1999), 9-74.

\bibitem{SW03} A. Soffer, M.I. Weinstein,
Selection of the ground state
for nonlinear Schr\"odinger equations,
{\em Rev. Math. Phys.} {\bf 16} (2004), no. 8, 977-1071.
arXiv:nlin/0308020.

\bibitem{Sp04} H. Spohn,
Dynamics of charged particles and their radiation field,
Cambridge University Press, Cambridge, 2004.

\bibitem{T02}
Tai-Peng Tsai,
Asymptotic dynamics of nonlinear Schr\"odinger
equations with many bound states,
{\em J. Differ. Equations} {\bf 192} (2003), no. 1, 225-282.








\end{thebibliography}
\end{document}